\documentclass[a4paper]{article}
\usepackage{graphicx}
\usepackage{amsmath,amsthm,amsfonts,amssymb}
\usepackage{newlfont}
\usepackage{amsfonts}
\usepackage{amstext}
\usepackage{color}
\usepackage{epsfig}

\theoremstyle{plain}
\newtheorem{thm}{Theorem}[section]
\newtheorem{lem}[thm]{Lemma}
\newtheorem{prop}[thm]{Proposition}
\newtheorem{defn}[thm]{Definition}
\newtheorem{rem}[thm]{Remark}
\newcommand{\R}{ \mathbb{R} }
\newcommand{\norm}[1]{\left\Vert#1\right\Vert}

\newcommand{\pf}{ \noindent {\bf Proof.} }

\begin{document}
\title{Optimal treatment planning governed by kinetic equations }%

\author{M. Frank\thanks{TU Kaiserslautern, PO Box 3049, D-67653 Kaiserslautern, GERMANY.} \and M.Herty \thanks{RWTH Aachen University, Templergraben 55, D-52056 Aachen, GERMANY.} \and A. N. Sandjo\thanks{
RWTH Aachen University, Templergraben 55, D-52056 Aachen, GERMANY.} }
 
\date{\today}
\maketitle
\begin{abstract}
In this paper we study a problem in radiotherapy treatment planning, where the evolution of the radiation field is governed by a deterministic Boltzmann transport equation. We show existence, uniqueness and regularity of solutions to an optimal dose distribution problem constrained by the Boltzmann Continuous Slowing-Down equation in an appropriate function space. The main new difficulty is the treatment of the stopping power term. Furthermore, we characterize optimal controls for problems governed by this transport equation.
\end{abstract}

\section{Introduction} 







%

Besides surgery and chemotherapy, the use of ionizing radiation is one of the main tools in the therapy of cancer today \cite{BelMai05,BelMai06,BelMai07,BelLiMai08}. According to WHO data, in the year 2007 there were about 11.3 million new cancer cases. More than half of the patients that are treated receive radiation therapy at one point during their treatment. Since the early days of radiation treatment high energy photons have been the most important type of radiation. Other types of radiation include high energy electrons and heavy charged particles like protons and ions. The latter type of radiation is of growing importance, but has not reached the widespread use of photons and electrons, yet. 
The aim of radiation treatment is to deposit enough energy in cancer cells so that they are destroyed. On the other hand, healthy tissue around the cancer cells should be harmed as little as possible. Furthermore, some regions at risk, like the spinal chord, should receive almost no radiation at all.

It is still current practice that treatment plans involve several fixed beam directions which are selected by an experienced physician by hand. Radiation facilities where the beam head rotates around the patient and where the beam is shaped by multileaf collimators are entering clinical practice. These methods have become known as Intensity-Modulated Radiation Therapy (IMRT). Patient motion during treatment is also one of the future challenges in the field of external beam radiotherapy. For instance, tumors near the lung move due to breathing. Techniques addressing this problem have become known as 4D radiotherapy (4DRT) \cite{BucBevRoa05}, meaning that time, the fourth dimension, also has to be taken into account. A further technique, named Image-Guided Radiotherapy (IGRT) is currently being developed. In this method, the radiation is used to create patient images during treatment. All of these novel techniques require mathematical modeling and optimization techniques \cite{GifHorWarFaiMou06,SheFerOliMac99}.

Before the treatment of the patient can be started, the expected dose distribution, i.e.\ the distribution of absorbed radiative energy in the patient, has to be calculated.  Most dose calculation algorithms in clinical use rely on the Fermi--Eyges theory of radiation. In recent work \cite{KriSau05}, however, it has been shown that these can produce errors of up to 12\% near inhomogeneities.

This work is based on dose calculation using a Boltzmann transport
equation \cite{HenIzaSie06,Lar97}. Similar to Monte Carlo simulations it relies
on a rigorous model of the physical interactions in human tissue
that can in principle be solved exactly. Monte Carlo simulations are
widely used, but it has been argued that a grid-based Boltzmann
solution should have the same computational complexity \cite{Bor98,Bor99}.
Furthermore, Monte Carlo can only be used in derivative--free methods
for optimal dose distributions. In constrast, when optimizing using Boltzmann's 
equations it is possible to exploit structural information for numerical and analytical 
purposes of the optimization problem \cite{TerKol02,FHS,HerPinSea07, TerKolVauHeiKai99}.  

Our starting point is the Boltzmann equation for particle transport in a medium:
consider a part of the patient's body which contains the region of
the cancer cells. We assume that this part of the body can be
described as a convex, open, bounded domain $ Z$ in $\mathbb{R}^3$.
Furthermore, we assume that $ Z$ has a smooth boundary with outward
normal vector $n$. The direction, into which the particle is moving
is given by $\Omega\in S^2$, where $S^2$ is the unit sphere in three
dimensions:
$$
\Omega\cdot\nabla_x \psi(x,\epsilon,\Omega) = \int_0^\infty \int_{S^2} \sigma_s(x,\epsilon',\epsilon,\Omega'\cdot\Omega)\psi(x,\epsilon',\Omega') d\Omega' d\epsilon' - \sigma_t(x,\epsilon) \psi(x,\epsilon,\Omega).
$$

Here, $\psi$ can be thought of as being the number of particles at $x \in\mathbb{R}^3$ with energy $\epsilon$, and direction $\Omega\in S^2$. Scattering is determined by the total cross section $\sigma_t$ and by the scattering kernel $\sigma_s$, which can be seen as the probabilty that a particle with initial energy $\epsilon'$ and initial direction $\Omega'$ has energy $\epsilon$ and direction $\Omega$ after the scattering event.

For high energy particles, small angle and energy changes are very likely, thus the scattering kernel $s$ is very forward-peaked. This fact is utilized to derive the Boltzmann Continuous Slowing-Down (BCSD) approximation \cite{LarMifFraBru97}. This model still allows large-angle scattering (which is important in radiotherapy applications) but describes energy-loss differentially.

\begin{equation}\label{EQ-BCSD}
\begin{split}
-\frac{\partial}{\partial\epsilon}(S(x,\epsilon)\psi(x,\epsilon,\Omega)) 
&+ \Omega\cdot\nabla_x \psi(x,\epsilon,\Omega)   \\
&= \int_{S^2} \sigma_s(x,\epsilon,\Omega'\cdot\Omega)\psi(x,\epsilon,\Omega') d\Omega' - \sigma_t(x,\epsilon) \psi(x,\epsilon,\Omega),
\end{split}
\end{equation}
where
$$
\sigma_s(x,\epsilon,\mu) = \int_0^\infty \sigma_s(x,\epsilon,\epsilon',\mu)d\epsilon'
$$
and the stopping power is
$$
S(x,\epsilon) = 2\pi \int_0^\infty \int_{-1}^1 (\epsilon-\epsilon')\sigma_s(x,\epsilon,\epsilon',\mu) d\mu d\epsilon'.
$$
This equation can be viewed as an initial boundary-value problem for \\
$(x,\epsilon,\Omega)\in Z\times (0,\infty)\times S^2$, $S^2$ being the unit sphere in $\R^3$. 
To formulate boundary conditions, we define the in- and outgoing boundaries as
\begin{equation*}
\Gamma_{\pm} := \{ (x,\Omega) \in  \partial Z\times S^2: n(x) \cdot \Omega >(<) 0 \}
\end{equation*}
and prescribe
$$
\psi(x,\Omega) = q(x,\Omega) \;\mbox{ on } \Gamma_{-}.
$$
The ``initial condition'' is 
$$
\psi(x,\infty,\Omega) = 0,
$$
meaning that there are no particles with arbitrary large energy. 

Different other approaches exist.  For a review on neutral particle codes that have 
been applied to the dose calculation problem we refer the reader to \cite{GifHorWarFaiMou06}.

A  number of functionals and methods have been devised to
describe the effect of radiation on biological tissue, cf.\ the
extensive lists of references in the reviews \cite{Bor99} and
\cite{SheFerOliMac99}. It is clear that the amount of destroyed
cells in a small volume, be they cancer or healthy cells, is not
directly proportional to the dose
\begin{equation}
  D(x) = \int_0^\infty \int_{S^2} S(x,\epsilon)\psi(x,\epsilon,\omega) d\omega d\epsilon
\end{equation}
deposited in that volume. However, no single accepted type of model
has emerged yet. Moreover, current biological models require input
parameters which are not known exactly~\cite{SheFerOliMac99}. This
is why the authors of ~\cite{SheFerOliMac99} opted not to
investigate these models but rather to focus on some general
mathematical cost functionals. A quadratic objective function
together with nonlinear constraints was identified as the most
versatile model.
Divide the domain into tumour tissue, normal tissue and a region at risk: $Z=Z_T\cup Z_N\cup Z_R$.
We prescribe a desired dose distribution $\bar D$, which usually has a constant value in $Z_T$ and is zero elsewhere. The problem of optimal treatment planning is to find an external beam distribution $q$ such that
$$
J = \frac{\alpha_T}{2}\int_{Z_T} (D-\bar D)^2 dx + \frac{\alpha_N}{2}\int_{Z_N} (D-\bar D)^2 dx + \frac{\alpha_R}{2}\int_{Z_R} (D-\bar D)^2 dx
$$
is minimial. Additionally, we might add the constraints
$$
D \geq D_\text{min} \quad\text{in $Z_T$ and}\quad D \leq D_\text{max} \quad\text{in $Z_R$},
$$
which ensure that all of the tumour tissue is sufficiently irradiatied and the region at risk receives a limited dose.

\section{ Main result }\label{sec:main result}

\newcommand{\X}{Z\times S^2}
\newcommand{\T}{ [ 0,\epsilon_{\max} ] }

The target area is modelled by a bounded, convex domain $Z \subset \R^3$ with smooth boundary 
$\partial Z.$ All results extend to the case $Z\subset \R^n$ and $S^{n-1}.$ 

We make some additional assumptions in order to provided a concise mathematical treatment. The first assumption is no major restriction. We assume that there is a maximal energy denoted by $\epsilon_{\max}$. Later, we will introduce a transformation $\epsilon' = \epsilon_\text{max}-\epsilon$ and solve an initial-boundary-value problem with initial values prescribed for $\epsilon'=0$.

Our second assumption is that the averaged scattering coefficients $\sigma_t$ and $\sigma_s$ do not depend on energy. This means that the elastic part of the scattering process is independent of the energy of the incident particles. This is not satisfied for elecron scattering. Our main purpose here is to deal with the stopping power. Thus, we postpone the treatment of energy-dependent scattering coefficients to future work. The stopping power itself is assumed to be independent of space. A possible reasoning for this assumption is that dose calculations are performed with data coming from a voxel--based CT scan. The stopping power is constant in each voxel and we couple the solutions to the transport equations over the different voxels.

Under these two assumptions we derive rigorous results on existence, uniqueness and regularity on optimal controls. We are thus interested 
in solutions to optimal control problems on $\X\times\T$  subject to the Boltzmann continuous slowing-down equation
in the following form
\begin{subequations}\label{bcsdmain} 
\begin{eqnarray}
 \partial_{\epsilon} S(\epsilon) \psi + \Omega \nabla_{x} \psi(x,\epsilon,\Omega) + \sigma_{t}(x,\Omega)
  \psi(x,\epsilon,\Omega) =  \\
  \int_{S^2} \sigma_{s}(x,\Omega' \cdot \Omega) \psi(x,\epsilon,\Omega') d\Omega' + q(x,\epsilon,\Omega), \\
 \psi(x,\Omega,\epsilon) = 0 \mbox{ on } \Gamma_{-}, \; \psi(x,\Omega,0) = 0 \mbox{ on } Z \times S^2.
\end{eqnarray}
\end{subequations}


\noindent We denote $x^+=\max(x,0)$ and we introduce the function spaces 
\begin{eqnarray*}
L^2_{ad} &=& \Big\{  q \in L^2( \X\times\T ): q \geq 0 \mbox{ a.e. } \Bigr\}, \\
D(A) &=& \{ \varphi \in L^2(\X): - \Omega\nabla_{x} \varphi \in L^2(\X), \; \varphi|_{\Gamma_{-}} = 0 \}, 
\end{eqnarray*}
and
\begin{subequations}\label{W2}
\begin{eqnarray*}
\mathbb{W}_{2} = \Big\{\psi\in L^2(\X\times\T));
    \frac{\partial\psi}{\partial \epsilon} +\Omega\cdot\nabla_{x}\psi \in L^2(\X\times\T);\\
     \psi(\cdot,\cdot,0) \in L^2(Z\times S^2);
      \psi|_{{\Gamma_{-}}\times \T}\in L^2(\Gamma_{-} \times \T,\Omega\cdot n(x) d\Gamma_{-})\Big\}.
\end{eqnarray*}
\end{subequations}
We impose the following assumptions on the coefficients: 
\begin{enumerate}
\item[{\bf A1}] The functions $\sigma_{t}$ and $\sigma_{s}$ are  non--negative a.e..
\item[{\bf A2}] The functions satisfy $\sigma_{t} \in L^{\infty}(Z)$ and $\sigma_{s} \in L^{\infty}( Z \times [-1,1] )$
\item[{\bf A3}] The scattering kernel is uniformly bounded for all $x \in Z:$
$$ \int_{-1}^1 \sigma_{s}( x,  \mu) d\mu \leq c.$$
\item[{\bf A4}] The stopping $S=S(\epsilon)$ is a strictly positive and continuous  function.
\end{enumerate}

\begin{figure}
\centering\includegraphics[width=0.8\linewidth]{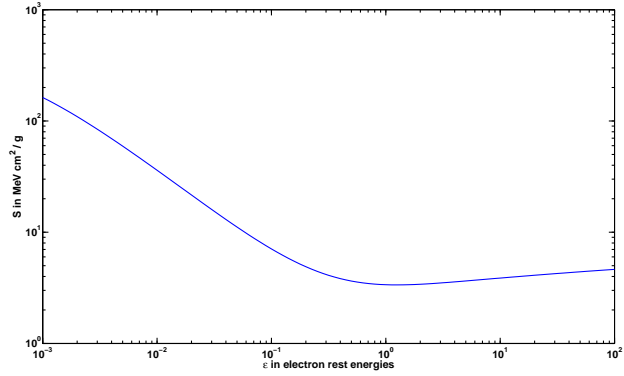}
\caption{Moller stopping power for water as a function of energy.}
\label{fig:SM}
\end{figure}

\begin{rem}
The assumptions (A1)--(A4) are satisfied for any reasonable scattering and absorption kernel. Figure~\ref{fig:SM} shows the stopping power for Moller inelastic scattering for water. In units of the electron rest energy, it can be written as \cite{HenIzaSie06}
\begin{equation*}
\begin{split}
S(\epsilon) = \frac{2\pi r_e^2\rho(\epsilon+1)^2}{\epsilon(\epsilon+1)}\Big( & \frac{\epsilon}{\epsilon-\epsilon_B} + 2\ln\frac{\epsilon-\epsilon_B}{2\epsilon_B(\epsilon-\epsilon_B)}\\ 
&+ \frac{1}{2(\epsilon+1)^2}( \frac{(\epsilon-\epsilon_B)^2}{4}-\epsilon_B^2 ) -\frac{2\epsilon+1}{(\epsilon+1)^2}\ln 2  \Big).
\end{split}
\end{equation*}
Here, $\epsilon_B$ is the electron binding energy in a water molecule.
\end{rem}

The positivity and smootheness of $S$ alow us to introduce yet another new variable $\tilde\epsilon = r(\epsilon)$ as the unique solution to 
$$
\frac{dr}{d\epsilon} = \frac{1}{S(\epsilon)},\quad r(0) = 0.
$$
Then the quantity $\tilde\psi(x,\tilde\epsilon,\Omega):= \psi(x,r^{-1}(\tilde\epsilon),\Omega)$ satisfies a transformed transport equation without stopping power term, but with a modified right hand side. By investigating this transformed equation, we obtain the following main result:
\begin{thm}\label{main} Assume (A1)--(A4). 
Let $\alpha_{i} \in L^{\infty}(Z)$ be positive a.e., let $\bar q \in L^2_{ad}$ and let  $\bar \psi \in  L^2( \X\times\T )$.  
\par
Then, the problem to minimize
$$ 
\int_{\T\times Z}  \frac{\alpha_{1}(x)}2  \left( \int_{S^2} \psi - \bar\psi d\Omega \right)^2 dxd\epsilon 
  + \int_{\T\times\X}\frac{\alpha_{2}(x)}2 (q - \bar q )^2 dxd\Omega d\epsilon
  $$
subject to 
\begin{eqnarray*} 
 \partial_{\epsilon} (S \psi) + \Omega \nabla \psi + \sigma_{t} \psi = \int_{S^2} \sigma_{s} \psi d\Omega' + q, \\
 \psi(x,\Omega,\epsilon) = 0 \mbox{ on } \Gamma_{-}, \; \psi(x,\Omega,0) = 0 \mbox{ on } Z \times S^2,
\end{eqnarray*}
admits a unique weak solution  $(\psi^*, q^*) \in C^0(\T,L^2(\X))\times L^2_{ad}$.  Under the assumption $\bar \psi \in   \mathbb{W}_{2}$ and  
$\psi^*\in \mathbb{W}_{2} \cap C^0(\T, D(A))$, $(\psi^*, q^*)$ is a local minimum, if there exists a weak solution $\lambda^* \in C^0(\T,D(A^*)) \cap \mathbb{W}_{2}$ of the first--order optimality system:
\begin{subequations}\label{opti main}
\begin{eqnarray}
\partial_{\epsilon} S \psi + \Omega \nabla \psi + \sigma_{t} \psi = \int_{S^2} \sigma_{s} \psi d\Omega + q, \\
 \psi(x,\Omega,\epsilon) = 0 \mbox{ on } \Gamma_{-}, \; \psi(x,\Omega,0) = 0 \mbox{ on } \X, \\
- S \partial_{\epsilon} \lambda - \Omega \nabla \lambda + \sigma_{t} \lambda = \int_{S^2} \sigma_{s} \lambda d\Omega' + \alpha_{1} \int_{S^2}\left( \psi - \bar \psi \right) d\Omega, \\
 \lambda(x,\Omega,\epsilon) = 0 \mbox{ on } \Gamma_{+}, \; \lambda(x,\Omega,\epsilon_{\max}) = 0 \mbox{ on } 
 \X, \\
q^* = \left( q^* - \lambda^*  - \alpha_{2} ( q^*-\bar q) \right)^+.
 \end{eqnarray}
 \end{subequations}
\end{thm}

\section{Existence and uniqueness of the minimizer}Ê
The proof of Theorem \ref{main} is split among several results. We start with the notion of a weak solution.
\begin{defn}\label{def:weak} 
Let $\mathcal{T}$ be a strongly continuous semigroup of operators on  a Banach space $\mathcal{X}$ with infinitesimal generator $\mathcal{A}$. Let $f\in L^1( [0,T], \mathcal{X})$, $T>0$,   $S:\R\to\R^+$ be a continuous, positive real--valued function and let $r$ be the solution to $\frac{dr}{dt} = \frac{1}{S(t)}$ with $r(0)=0.$ 

We call $\varphi \in C^0( [0,T], \mathcal{X})$ a weak solution to 
\begin{equation}\label{general:def} \partial_{t} \varphi 
= \mathcal{A} \varphi  + f(t), \varphi(0) = \varphi_{o},\end{equation} iff $\varphi (t) = \mathcal{T}(t) \varphi_{o} +\int_{0}^t 
\mathcal{T}(t-s)f(s) ds.$

We call  $\phi \in C^0( [0,T], \mathcal{X})$  a weak solution to $$\partial_{t}S(t) \phi 
= \mathcal{A} \phi  + \tilde f (t), \phi(0) = \phi_{o}, $$ iff $\varphi( t ) = S( r^{-1}(t)  ) \phi( r^{-1}(t) )$ is a weak solution to (\ref{general:def})
with $f( t ) =\tilde f( r^{-1}(t) ) S( r^{-1}(t) )$ and $\varphi_{o}=S(0)\phi_{o}$. We call $\omega \in C^0( [0,T], \mathcal{X})$ a weak solution to 
$$ S(t) \partial_{t} \omega = A\omega + \tilde f(t), \omega(0)=\omega_{o}, $$ iff $\varphi(t) = \omega( r^{-1}(t) ), \varphi_{o}=\omega_{o}$ is a weak solution to (\ref{general:def}). 
\end{defn}
\begin{rem}\label{rem}
The given definition (\ref{general:def}) is as the definition of weak solutions in \cite[Volume 5]{RL} or, equivalently, to the definition of  mild solutions  as given e.g. in~\cite{HA}.  Since $S$ is positive, the function $r$ is invertible and due to the initial condition $r(0)=0$.
In the case $S\equiv 1$ all definitions coincide.  The motivation for the latter definitions is as follows:  Provided we have sufficient regularity we can compute $\partial_{t} \left(  S(t)\phi(t) \right)  = \frac{ \partial_{r} \varphi( r(t) )}{ S(t) } = \mathcal{A} \phi + \tilde f(t)$ and, similarly, we compute $S(t) \partial_{t} \omega(t) = \partial_{r} \varphi ( r(t) )= \mathcal{A} \omega + f(t).$ Similarly,  a weak solution on $[0,T_{\max}]Ê\times \mathcal{X}$
$$ - \partial_{t} \tilde \phi = \mathcal{A}\tilde \phi + \tilde f(t), \; \tilde \phi(T) = \phi_{o}$$ 
is defined by a weak solution $\phi(t) = \phi(T_{\max}-t)$ to (\ref{general:def}) with $f(t)= \tilde f(T-t).$
\end{rem}
Next, we state semigroup properties for (\ref{bcsdmain}) in the case $S\equiv 1.$ Denoting by $\mathcal{L}(\mathcal{X},\mathcal{Y})$ the set of linear and bounded operators,  under the assumptions (A1)--(A3), we have~\cite{RL,FHS}
$$ (\Sigma \varphi) := \sigma_t(x,\Omega)\varphi \in  \mathcal{L} \left( L^2(\X) , L^2(\X) \right),$$ 
$$ (K\varphi) = \int_{S^2}\sigma_s(x,\Omega'\cdot\Omega)\varphi(x,\Omega' )d\Omega' \in \mathcal{L}\left( L^2(\X), L^2(\X) \right)$$ and 
$$ (A\varphi) = - \Omega \nabla_{x} \varphi : D(A) \to L^2(\X)$$ 
is an unbounded operator on $L^2(\X)$ with domain of definition $D(A)$. We define
\begin{equation}\label{defT}ÊT := A-\Sigma+ K : D(A)\subset L^2(\X) \to L^2(\X)\end{equation} 
and the general problem  (\ref{bcsdmain}) reads 
\begin{equation}\label{bcsd-2} 
\partial_{\epsilon} S \psi = T \psi + q, \;  \psi  = \psi_{0} \in \X\times\{0\} \mbox{ and } \psi = 0 \in \Gamma_{-}\times\T .
\end{equation}
We summarize the properties of the defined operators in Proposition~\ref{p1}, see~Theorem~2 in Volume 6, Chap. XXI~\cite{RL} and Definition~1 and Theorem~1 in Volume 6, Chap. XXI for definition and properties of the trace of $\psi \in D(A).$ Due to Remark~3 in Volume 6, Chap. XXI,  we have that  $\psi \in D(A)$ implies $ \psi \in L^2( \Gamma_{-} |\Omega \cdot n(x) | d\Gamma_{-} ),$ where $n(x)$ is the outer normal at $x\in \partial Z.$Ê  
\begin{prop}\label{p1}
Assume (A1)--(A3). Then, 
\begin{enumerate}
\item[(1)] $A$ is the infinitesimal generator of a strongly continuous semigroup in $L^2(\X).$ 
\item[(2)] $T$ is the infinitesimal generator of a strongly continuous semigroup in $L^2(\X)$ and $D(T) = D(A).$
\item[(3)]  For any $\psi_0 \in D(T)$ there exists a unique classical solution $\psi \in C^{1}(\T;L^2(\X))\cap C^0(\T;D(T))$ 
of the Cauchy problem 
$$ \partial_{\epsilon} \psi = T \psi, \; \psi(0) = \psi_0, \; \psi = 0 \mbox{ on }Ê\Gamma_{-}\times \T$$ 
The solution is given by 
$\psi(\epsilon)= \exp(\epsilon T)\psi_0$ where $\exp(\epsilon T)$ is the strongly continuous semigroup with generator $T$. 
\end{enumerate}
\end{prop}
For the non--homogenous problem (\ref{bcsd-2}) and $S\equiv 1$ the following result is classical~\cite{HA,RL} and in fact, it sufficies to verify the regularity properties of $\psi$ as defined below. 
\begin{prop}\label{p2}
Assume (A1)--(A3),  let $q\in L^2(\T\times\X)$ and $\psi_0\in L^2(\X).$ Then, for $S\equiv 1$, there exists a unique weak solution $\psi \in C^0( \T, L^2(\X))$ to   given by $\psi = \exp( \epsilon T)\psi_0 +\int_{0}^{\epsilon} \exp( (\epsilon - s)T ) q(s) ds.$ 

\par

If additionally $\psi_0\in D(T)$ and $q\in C^1(\T,L^2(\X)),$ then  $\psi$ is a classical solution to (\ref{bcsd-2}). We have $\psi \in C^1(\T,L^2(\X)) \cap C^0(\T,D(T)).$

If $q\geq 0$ and $\psi_0 \geq 0,$ then we obain $\psi \geq 0.$
\end{prop}
\pf 
A proof of the second and third statement is given by  Theorem~3 in Volume 6, Chap. XXI of \cite{RL}.  The first statement is classical and a proof is found Lemma 11.14, Theorem 11.16 in \cite{RR}. $\hfill \blacksquare$\\
Under the given assumptions on $\psi_0$ and $q$ we do not  necessarily obtain $\psi \in C^0(\T, D(T)).$ However,
if $\psi_0\in D(T)$ and $q\in L^2(\X),$ then  $\psi \in C^1(\T,L^2(\X)) \cap C^0(\T,D(T)).$ Additionally, we 
have (Theorem~3, Chap. XXI, Volume 6, \cite{RL}), that if $\psi$ belongs to $\mathbb{W}_{2}$
then it is a pointwise a.e. solution to (\ref{bcsdmain}) and it is unique in this space: 
For functions $\psi \in \mathbb{W}_{2}$ we may apply Green's formula to obtain for any $\tau \in \T:$
\begin{eqnarray*}
 2\int_0^{\tau}{\left\langle \psi(\epsilon),K\psi(\epsilon)-\Sigma\psi(\epsilon)\right\rangle}d\epsilon=\\
{\norm{\psi(\tau)}^2}_{L^2(\X)}+\int_0^{\tau}\left({\norm{\psi(\epsilon)}^2}_{L^2(\Gamma_{+})}-{\norm{\psi(\epsilon)}^2}_{L^2(\Gamma_{-})}\right)d\epsilon\geq{\norm{\psi(\tau)}^2}_{L^2(\X) }
\end{eqnarray*}
and hence  ${\norm{\psi(\tau)}^2}_{L^2(Z\times S^2)}\leq 2 \norm{K}\int_0^{\tau}{\norm{\psi(\epsilon)}^2}_{L^2(Z\times S^2)}d\epsilon$
and yields uniqueness due to Gronwall's lemma. The previous Proposition~\ref{p2} also allows to define a control--to--state operator at first in the case $S\equiv 1$ by 
\begin{equation}\label{Xi} \Xi (q; \psi_0 ) = \psi  = \exp( \epsilon T)\psi_0 +\int_{0}^{\epsilon} \exp( (\epsilon - s)T ) q(s) ds. \end{equation}

\begin{lem}\label{lem1}
Assume (A1)--(A3) and let $q \in L^2(\X\times\T)$ and $\psi_0 \equiv 0.$ Then, the operator $\Xi$ is a linear and bounded operator
from $L^2(\X\times\T) \to L^2(\X\times\T).$
\end{lem}
\pf Note that for any strongly continuous semigroup $\exp( \epsilon T) $ there exists a constant $\omega\geq 0$ and $M\geq 1,$ such that 
$\| \exp(\epsilon T)\|\leq M \exp(\epsilon \omega ),$ see \cite{PA}. $\Xi$ is linear and since $C^0(\T,L^2(\X))\subset L^2(\T\times\X)$ we have 
$\norm{\Xi(q)}_{L^2(\X\times\T} \leq Me^{\omega \epsilon_{\max}}\cdot{\norm q}_{L^2(\X\times\T)}.$           $\hfill\blacksquare$\par             
This yields existence of a minimizer for a class of cost functionals. We obtain the following result as extension to Theorem 3.1 in \cite{FHS}.
\begin{thm}\label{existence1}
Assume (A1)--(A3). Let $\alpha_{i} \in L^{\infty}(Z)$ be positive a.e., let $\bar \psi\in L^2(\X\times\T)$ and $\bar q\in L^2_{ad}.$ Then, the problem 
\begin{subequations}\label{opti system}
\begin{eqnarray}
&& \min \int_{\T\times Z} \frac{ \alpha_{1} }2 \left( \int_{S^2} \psi - \bar \psi \right)^2 dxd\epsilon + \int_{\T\times\X} \frac{\alpha_{2}}2 ( q- \bar q)^2 dxd\Omega d\epsilon  \\ 
&& \mbox{ subject to } \\
&& \partial_{t} \psi + \Omega \nabla_{x} \psi + \sigma_{t} \psi = \int_{S^2} \sigma_{s} \psi d\Omega' + q, \\
&& \psi = 0 \mbox {on }Ê\Gamma_{-},  \; \psi(x,\Omega,0) = 0 \mbox{ on } \X
\end{eqnarray}
\end{subequations}
admits a unique minimizer $q^* \in L^2_{ad}$ and corresponding weak solution $\psi^* \in C^0(\T,L^2(\X)).$ 
\end{thm}
\pf 
For given $\bar\psi \in L^2(\X)$ the operator 
\begin{equation}\label{phi}
\Phi:L^2(\X\times\T) \to L^2(Z\times\T)
\end{equation}
 given by $\Phi(\psi) = \int_{S^2} \psi -\bar \psi d\Omega$ is an affine linear, bounded operator. Due to Lemma~\ref{lem1} the operator $\Xi$ is a linear and bounded operator  from $L^2_{ad} \subset L^2(\X\times\T) \to L^2(\X\times\T).$ The subspace $L^2_{ad}$ is a  closed subset of a Hilbert space. Due to Theorem 2.16 in \cite{FT}  the convex functional  
$$\int_{\T\times Z} \frac{ \alpha_{1} }2 \Phi (\Xi(q) )^2 dxd\epsilon + \int_{\T\times\X} \frac{\alpha_{2}}2 ( q- \bar q)^2 dxd\Omega d\epsilon$$
attains its minimum. The latter being unique provided that $\alpha_{2}>0.$ $\hfill\blacksquare$\par  

Similarly, we obtain an existence result for the equation with $S\not = 1.$
\begin{thm}\label{existence2}
Assume (A1)--(A4). Let $\alpha_{i} \in L^{\infty}(Z)$ be positive a.e., let $\bar \psi\in L^2(\X\times\T)$ and $\bar q\in L^2_{ad}.$ Then, the problem 
\begin{eqnarray*}
\min \int_{\T\times Z} \frac{ \alpha_{1} }2 \left( \int_{S^2} \psi - \bar \psi \right)^2 dxd\epsilon + \int_{\T\times\X} \frac{\alpha_{2}}2 ( q- \bar q)^2 dxd\Omega d\epsilon  \\ \mbox{ subject to } \\
\partial_{t} S \psi + \Omega \nabla_{x} \psi + \sigma_{t} \psi = \int_{S^2} \sigma_{s} \psi d\Omega' + q, \\
\psi = 0 \mbox {on }Ê\Gamma_{-},  \; \psi(x,\Omega,0) = 0 \mbox{ on } \X
\end{eqnarray*}
admits a unique minimizer $q^* \in L^2_{ad}$ and corresponding weak solution $\psi^* \in C^0(\T,L^2(\X)).$ 
\end{thm}
\pf   Due to (A4), the solution $\frac{d}{dt} r= \frac{1}{S(t)}, r(0)=0$ is a strictly monotone, smooth function with smooth inverse $r^{-1}$. 
For any given $q(x,\Omega,\epsilon) \in L^2_{ad}$ introduce the operator $\tilde \; : L^2_{ad}\to L^2_{ad}$ by 
$$\tilde q(x,\Omega,\epsilon) := q(  x, \Omega, r^{-1}(\epsilon) ) S( r^{-1} (\epsilon) )$$ and we  denote by $\phi = \Xi( \tilde q )$ for $\Xi$ defined in (\ref{Xi}).  Due to Proposition \ref{p2}  $\phi\in C^0(\T,L^2(\X))$ exists. Due to Definition~\ref{def:weak}  $\psi(x,\Omega,\epsilon) = \phi(  x, \Omega, r(\epsilon) ) / S (\epsilon )$ is a weak solution to equation (\ref{bcsdmain}) with zero boundary conditions since $r(0)=0.$
 We have $\psi \in C^0(\T,L^2(\X))$ and we define the operator   
\begin{equation} \label{X} X :L^2_{ad} \to L^2(\X\times\T): X(q)(x,\Omega,\epsilon) =  \frac{\Xi\left( \tilde q ) \right) ( x, \Omega, r(\epsilon)) }{S(\epsilon)}.\end{equation}
Due to Lemma~\ref{lem1}, the continuity of  $S$ on the closed set $\T$ and the $C^1-$property of $r$, the operator $X$ is a linear bounded operator. Then, the proof is exactly as in the Theorem \ref{existence1} when replacing $\Xi$ by $X.$  $\hfill\blacksquare$Ê\par

\section{First-order optimality conditions}
We define the operator 
\begin{equation}\label{adjT} T^*:= - A - \Sigma + K:  D(A^*) \subset L^2(\X) \to L^2(\X) \end{equation} 
for
$$ D(A^*) := \{ \lambda \in L^2(\X):  \Omega \nabla_{x} \lambda \in L^2(\X), \;  \lambda|_{\Gamma^+} = 0 \}$$  
and study for some $r\in L^2(\X\times\T)$ the equation
\begin{equation}\label{bcsd adjoint}
 - \partial_{\epsilon} \lambda = T^* \lambda + r, \; \lambda(x,\Omega,\epsilon_{\max}) = 0, \; \lambda = 0 \mbox{ on }Ê\Gamma^+. 
\end{equation}
\begin{lem}\label{lem:adj}
Under the assumptios (A1)--(A3), there exists a linear and bounded operator $\Xi^*:L^2(\X\times\T) \to L^2(\X\times\T)$ with $\lambda=\Xi^*(r)$ being the weak solution to (\ref{bcsd adjoint}) for any $r\in L^2(\X\times\T)$. Additionally, we have $\lambda \in C^0(\T, L^2(\X)).$  The operator $\Xi^*$ is the adjoint operator on $L^2(\X\times\T)$ to $\Xi$ provided that  $\Xi, \Xi^* :  L^2(\X\times\T) \to \mathbb{W}_{2} \cap C^0(\T, D(A^{(*)})) \subset L^2(\X\times\T).$ 
\end{lem}
\pf Under the given assumptions and due to Theorem~XX in Volume 6, Chap. XXI \cite{RL}, the operator $T^*$ is the infinitesimal generator of a strongly continuous semigroup  with domain of definition $D(A^*).$ Hence, the equation $\partial_{\epsilon} \mu = TÊ\mu + f$  with zero initial data  $\mu(0)=0$ admits a weak solution given by  
$\mu = \int_{0}^{\epsilon} \exp( (\epsilon-s)T^* )f(s) ds \in C^0(\T,L^2(\X))$ for any $f\in L^2(\X\times\T).$  Hence, $\lambda(\epsilon) = \mu( \epsilon_{\max}  - \epsilon) \in C^0(\T,L^2(\X)) $ is a weak solution to (\ref{bcsd adjoint}).  The solution $\lambda$ can be written as
$\lambda(\epsilon) = - \int_{\epsilon_{\max}-\epsilon}^{\epsilon_{\max}} \exp( ( \epsilon_{\max}-\epsilon - s)T^*) f(s) ds.$
This defines a solution operator $\Xi^*$ and as in Lemma~\ref{lem1} the operator $\Xi^*$ is linear and bounded. Due to Theorem 3.3 in \cite{FHS}
the operator $T^*$ is adjoint to $T$ on $L^2(\X\times\T).$ Given any $ z,w \in L^2(\X\times\T)$, we denote by  $ \lambda=\Xi^{*}(z) \in \mathbb{W}_{2}$  and 
by $\psi=\Xi(w) \in \mathbb{W}_{2}$.  Note that for $\lambda, \psi \in \mathbb{W}_{2}$ we may apply Green's formula and hence obtain 
\begin{eqnarray*}
\begin{split}
<w, \Xi^*(z)>_{L^2(\X\times\T)} = \int_0^{\epsilon_{\max}} {\left[{\left\langle \psi(\epsilon),-\Omega\cdot\nabla_{x}\lambda(\epsilon)+\frac{\partial \lambda(\epsilon)}{\partial \epsilon}\right\rangle}_{L^2(\X)}\right]}d\epsilon + \\
\int_0^{\epsilon_{\max}} {\left[{\left\langle K\psi(\epsilon)-\Sigma\psi(\epsilon),\lambda(\epsilon)\right\rangle}_{L^2(\X)}\right]}d\epsilon 
=\\ 
\left\langle \psi(\epsilon_{\max}),\lambda(\epsilon_{\max})\right\rangle_{L^2(\X)}-\left\langle \psi(0),\lambda(0)\right\rangle_{L^2(\X)}+\\
\int_0^{\epsilon_{\max}} \left[{\left\langle\psi(\epsilon),\lambda(\epsilon)\right\rangle_{\Gamma_{+}}-\left\langle\psi(\epsilon),\lambda(\epsilon)\right\rangle_{\Gamma_{-}}}\right]d\epsilon = <\Xi(w), z >_{L^2(\X\times\T)}.
\end{split}
\end{eqnarray*} $\hfill\blacksquare$Ê\par
The previous result can be used to deduce the existence of a first--order optimality system in the case $S\equiv 1.$
This result will be extended to the more general case below. 
\begin{thm}\label{opti1}
Assume (A1)--(A3) and $\bar\psi \in \mathbb{W}_{2}.$ Further,  assume that the minimizer $(q^*, \psi^*=\Xi(q^*) )$ of (\ref{opti system}) belongs to $L^2_{ad} \times \mathbb{W}_{2}.$  Then, the first--order necessary optimality conditions are  $$q^* = ( q^* - \lambda^* - \alpha_{2} ( q^* - \bar q)  )^+ \; \mbox{ a. e. } \in \X\times\T,$$  
provided that  $\lambda^* = \Xi^*( \alpha_{1} \int_{S^2} \psi^* - \bar \psi d\Omega)$ and $\lambda^*Ê\in \mathbb{W}_{2} \cap C^0(\T, D(A^*))$.
\end{thm}
\pf  The proof is similar to the proof of Theorem 3.3 \cite{FHS}. Consider the case $\bar \psi \equiv 0$ first. Then, 
the problem (\ref{opti system}) can be rewritten as 
$$Ê\min \| \sqrt\frac{\alpha_{1}}{2} \Phi( \Xi ( q ) ) \|_{L^2(Z\times\T)}^2 + \| \sqrt{ \frac{\alpha_{2}}2 } (q-\bar q) \|_{
L^2(\X\times \T )}^2 =: j(q) $$
where $\Phi$ is defined in (\ref{phi}) and $\Xi$ defined in (\ref{Xi}).  The adjoint operator on $L^2(Z\times\T)$ to $\Phi:L^2(\X\times\T) \to L^2(Z\times\T)$ satisfies $< w, \Phi(v) >_{L^2(Z\times \T)}  = <\Phi^*(w), v>_{L^2(\X\times\T)} =  \int_{\X\times\T} w\; v \; dxd\Omega d\epsilon$. 
Due to Lemma 2.20, Lemma 2.21 and Theorem 2.22 of \cite{FT} 
and Theorem 3.3 in \cite{FHS}, we obtain the gradient by  $ j'(q)  = \left(   \Phi( \Xi )
\right)^{*} \left( \alpha_{1} \Phi(\Xi (q) ) \right)  + \alpha_{2} ( q -\bar q) .$ Under the given regularity assumptions  and  denoting by $\lambda^* = \Xi^*( \alpha_{1} ( \psi^*  )$ and $\psi^*=\Xi(q^*)$ we rewrite the gradient as 
$$j'(q) =  \Xi^* ( \alpha_{1} \int_{S^2} \psi d\Omega) + \alpha_{2} ( q- \bar q) =  \lambda  + \alpha_{2}( q- \bar q).$$ In the case $\bar \psi \not = 0$ we consider $\varphi$ such that $\Xi( \varphi ) = \bar \psi$ and as in Theorem 3.3 \cite{FHS} we consider an optimal control problem for $\tilde \psi = \psi - \bar \psi$ subject  to $\Xi(q-\varphi) = \tilde \psi.$   Since $\Xi$ is linear and together with the result in the case $\bar \psi \equiv 0$  we obtain  the assertion. The strong form of the necessary conditions is given by system (\ref{opti system}) in the case $S\equiv 1.$
$\hfill\blacksquare$Ê\par
\par \noindent
{\bf Proof of Theorem~\ref{main} }.
Existence of a minimizer is given by Theorem~\ref{existence2}. The minimizer $q^*\in L^2_{ad}$ is unique provided that $\alpha_{2}>0$ and $\psi^*$ is a weak solution in the sense of Definition~\ref{def:weak}. It is unique in the space $\mathbb{W}_{2}.$ \par

Let $r$ be the solution to $r'(\epsilon) = \frac{1}{S(\epsilon)}, r(0)=0$. Since $S>0$ and continuous, the function $r$ is $C^1$ on $\T$ and 
strictly monotone. Therefore, $r$ is invertible and  $r' >0.$  For $q\in L^2(\X\times\T),$ and $T$ defined by (\ref{defT}), let  $\psi$ be the weak solution to (\ref{bcsdmain}), i.e., 
\begin{equation}\label{temp1} 
\partial_{\epsilon} S \psi = T\psi  + q, \; \psi(x,\Omega,0) =0, \; \psi = 0 \mbox{ on } \Gamma_{-}.
\end{equation}
Under the given assumptions and for any $q\in L^2(\X\times\T)$ there exists a weak solution $\psi \in C^0(\T,L^2(\X)) \subset L^2(\X\times\T)$  unique in $\mathbb{W}_{2}$  due to Theorem~\ref{existence2}. Then, 
$$\phi(x,\Omega,r(\epsilon)) / S(\epsilon) :=  \psi(x,\Omega,\epsilon), \; \tilde q(x,\Omega,\epsilon) = q(x,\Omega, r^{-1}(\epsilon)) S(r^{-1}(\epsilon))$$
is a weak solution to 
\begin{equation}\label{temp2}
\partial_{\epsilon}  \phi = T\phi  + \tilde q, \; \phi(x,\Omega,0) =0, \; \phi = 0 \mbox{ on } \Gamma_{-}.
\end{equation}
The solution operator is denoted by $\Xi$ as in (\ref{Xi}). Similarly, for $(\phi,\tilde q)$ being a weak solution to (\ref{temp2}), there exists a unique $(\psi, q)$ being a weak solution to (\ref{temp1}). Hence, there exists a bijective mapping on $C^0(\T, L^2(\X)) \times L^2_{ad}$  from solutions$(\psi,q)$ to (\ref{temp2}) to weak solution $(\phi,\tilde q)$ to (\ref{temp1}). The mapping preserves the regularity of the weak solutions, i.e., if $\psi \in\mathbb{W}_{2} \cap C^0(\T, D(A))$ then $\phi \in \mathbb{W}_{2} \cap C^0(\T, D(A)).$ Hence,  we have for $T_{R} = r(\epsilon_{\max}).$ 
 \begin{eqnarray} 
\int_{\T\times Z}  \frac{ \alpha_{1} }2  \int_{S^2} \left( \psi(x,\Omega,\epsilon) - \bar \psi(x,\Omega,\epsilon) d\Omega \right)^2  + \nonumber 
\\ \label{temp3} \\
\nonumber
 \int_{\T\times\X} 
 \frac{ \alpha_{2} }2 \left( q - \bar q \right)^2(x,\Omega,\epsilon)  dx d\Omega d\epsilon =\\
 \nonumber
 \int_{\T\times Z}  \frac{ \alpha_{1} }2  \int_{S^2} \left( \frac{\phi(x,\Omega,r(\epsilon))  }{S(\epsilon)} - \bar \psi(x,\Omega,\epsilon) d\Omega \right)^2  + \\ \nonumber
 \nonumber \int_{\T\times\X} 
 \frac{ \alpha_{2} }2 \left( q - \bar q \right)^2(x,\Omega,\epsilon)  dx d\Omega d\epsilon = \\
\nonumber
 \int_{0}^{T_{R}} \int_{Z}Ê \frac{ \alpha_{1} }2  
\left( \int_{S^2} \frac{\phi(x,\Omega,\tau ) 
}{S(r^{-1}(\tau))}  - \bar \psi(x,\Omega,r^{-1}(\tau) ) d\Omega \right)^2  \frac{ d\tau dx}{ r'( r^{-1}(\tau)  ) }  + 
\\ \label{temp4}
\\ \nonumber
 \int_{0}^{T_{R}}\int_{\X}\frac{ \alpha_{2} }2 \left( \frac{ \tilde q(x,\Omega,\tau) }{ S(r^{-1}(\tau)) } - \bar q(x,\Omega,r^{-1}(\tau)) \right)^2  \frac{ d\tau dx d\Omega}{ r'( r^{-1}(\tau)  ) } 
\end{eqnarray}
where $\phi,\tilde q$ are weak solutions to (\ref{temp2}). Hence, if $(\psi,q)$ of (\ref{temp3}) yields  a minimizer $(\phi,\tilde q)$ for (\ref{temp4}). 
To obtain the necessary conditions, we may apply  Theorem~\ref{opti1} to  
$$  \alpha^0_{i}  := \frac{\alpha_{i}}{ S(r^{-1}(\tau)) },  \; \bar \psi^0:= S(r^{-1}(\tau)) \bar \psi(x,\Omega, r^{-1}(\tau), \bar q^0 = S(r^{-1}(\tau)) \bar q(x,\Omega, r^{-1}(\tau) ) $$   
since $\alpha^0_{i} \in L^{\infty}(\X\times \T)$ and $ \alpha^0_{i}> 0$ a.e., $\bar q^0 \in L^2(\X\times\T)$ and $\bar \psi^0 \in \mathbb{W}_{2}.$ We obtain on $[0,T_{R}] \times \X:$
\begin{eqnarray*}
q^0(x,\Omega,\tau) = \left( q^0(x,\Omega,\tau) - \lambda^0(x,\Omega,\tau) - \alpha_{2}  ( \frac{ q^0(x,\Omega,\tau) }{ S(r^{-1}(\tau)) } - \bar  q(x,\Omega, r^{-1}(\tau) ) \right)^+, \\
\lambda^0(x,\Omega,\tau)  = \Xi^*(  \alpha_{1}  \int_{S^2} \frac{ \phi^0(x,\Omega,\tau) }{S(r^{-1}(\tau))}  - \bar \psi(x,\Omega,r^{-1}(\tau) )  d\Omega ), \\
\phi^0(x,\Omega,\tau) = \Xi (  q^0(x,\Omega,\tau) ).
\end{eqnarray*}
Finally, we rewrite the optimality system in terms of $\psi$ and $q$. From the last equation we obtain that $q(x,\Omega, \epsilon) := q^0( x,\Omega, r(\epsilon)) / S ( \epsilon )$ and $\psi(x,\Omega,\epsilon):= \phi^0(x,\Omega,r(\epsilon)) / S(\epsilon)$ is a weak solution to (\ref{temp1}) on $\T\times\X.$  The second equation gives a weak solution to the following adjoint equation with zero terminal conditions  and zero boundary conditions on $\Gamma_{+}$
\begin{equation}\label{temp5}
 - \partial_{\tau} \lambda^0(x,\Omega,\tau) = T^*  \lambda^0(x,\Omega,\tau)  + \alpha_{1} \int_{S^2} \psi(x,\Omega, r^{-1}(\tau)) - \bar \psi (x,\Omega, r^{-1}(\tau)) d\Omega
\end{equation}
Due to Remark~\ref{rem} and Definition~\ref{def:weak} we conclude from (\ref{temp5}) that 
$ \lambda(x,\Omega,r^{-1}(\tau)) := \lambda^0(x,\Omega,\tau)$
is a weak solution on $\T$ to 
\begin{equation}\label{temp6}
- S(\epsilon) \partial_{\epsilon} \lambda(x,\Omega,\epsilon) = T^* \lambda(x,\Omega,\epsilon) + \alpha_{1} \int_{S^2} \psi(x,\Omega,\epsilon) - \bar \psi(x,\Omega,\epsilon) d\Omega.\end{equation} 
Further, $\lambda$ as the same regularity as $\lambda^0.$  Finally, we reformulate the first equation of the optimality system. The equation is equivalent  a.e. in $[0,T_{R}] \times \X$ to  
\begin{eqnarray*}
\lambda^0(x,\Omega,\tau) - \alpha_{2} \left( \frac{ q^0(x,\Omega,\tau)}{S(r^{-1}(\tau))} - \bar q(x,\Omega,(r^{-1}(\tau)) \right)  \;  \begin{pmatrix}Ê = 0 & \mbox{ iff } q^0(x,\Omega,\tau) \geq  0  \\  \geq 0 & \mbox{ iff }Êq^0(x,\Omega,\tau) = 0 \end{pmatrix}Ê
\end{eqnarray*}
Since $S>0$ and $r$ bijective,  we obtain using the functions $\psi$ and $\lambda$  for a.e. in $\T \times \X$ 
\begin{equation}\label{temp7}
  q(x,\Omega, \epsilon ) = \Bigl(  q(x,\Omega,\epsilon ) - \lambda(x,\Omega,\epsilon) - \alpha_{2}  \left( q(x,\Omega,\epsilon ) - \bar q(x,\Omega, \epsilon ) \right) \Bigr)^+
\end{equation}
The equations (\ref{temp1}, \ref{temp6}, \ref{temp7}) comprise the first--order optimality system provided that $\lambda, q$ and $\psi$ and therefore $\phi, \lambda^0$ and $q^0$  fulfill the  given regularity assumptions.  $\hfill\blacksquare$

\section{Remarks and further discussion}
We offer the following remarks and notes on further discussion to the results of Theorem \ref{main}. 
\begin{itemize}

\item 
The given results extend to minization problems of the type 
\begin{eqnarray}\label{cfct1}
\frac{\alpha_1}{2}\int_0^{\epsilon_{\max}}\int_{Z\times S^2}(\psi-\overline{\psi})^2dxd\Omega d\epsilon+\frac{\alpha_2}{2}\int_0^{\epsilon_{\max}}\int_{Z\times S^2}(q-\overline{q})^2dxd\Omega d\epsilon
\end{eqnarray}
for given functions $\overline{\psi},\overline{q} \in L^2(Z\times S^2\times(0,\epsilon_{\max}))$ by setting $\Phi \equiv Id.$ 
\par
Under additional assumptions on the minimizer $q^*$ the assertions on $\psi^*$ and $\lambda^*$ can be obtained, e.g., if $q^* \in W^{1,2}(\T, D(A))$ we obtain that $\Xi(q^*) \in \mathbb{W}_{2} \cap C^0(\T, D(A))$ and similarly for $\lambda^*.$

\item
Formally, we obtain the assertions of Theorem \ref{main} as follows. Let $\varphi(x,\Omega,r(\epsilon))= \psi(x,\Omega, \epsilon)$ with $r'(\epsilon) = \frac{1}{S(\epsilon)}$ and $r(0)=0.$ Then, the functional is 
\begin{eqnarray*}
J(\phi,q) =\int_{\T\times Z}  \frac{ \alpha_{1} }2  \int_{S^2} \left( \frac{\phi(x,\Omega,r(\epsilon))  }{S(\epsilon)} - \bar \psi(x,\Omega,\epsilon) d\Omega \right)^2  \\
+ \int_{\T\times\X} 
 \frac{ \alpha_{2} }2 \left( q - \bar q \right)^2(x,\Omega,\epsilon)  dx d\Omega d\epsilon.
\end{eqnarray*}
Using a coordinate transformation in $\epsilon$ by  $\tau = r(\epsilon)$ and  $T_{R} = r(\epsilon_{\max})$ we obtain  
\begin{eqnarray*}
J(\phi,q) =\int_{0}^{T_{R}} \int_{Z}Ê \frac{ \alpha_{1} }2  
\left( \int_{S^2} \frac{\phi(x,\Omega,\tau ) 
}{S(r^{-1}(\tau))}  - \bar \psi(x,\Omega,r^{-1}(\tau) ) d\Omega \right)^2  d\tau dx+ \\
 \int_{0}^{T_{R}}\int_{\X}\frac{ \alpha_{2} }2 \left( q - \bar q \right)^2 \left(x,\Omega,r^{-1}(\tau) 
 \right) \frac{ d\tau dx d\Omega}{ r'( r^{-1}(\tau)  ) } 
\end{eqnarray*}
and its formal derivative in direction $\delta\phi$ as 
\begin{eqnarray*} 
 \frac{d}{d\phi} J(\phi,q) \delta \phi =  \\
 \int_{0}^{T_{R}} \int_{\X} \alpha_{1} \left( \int_{S^2}  \frac{ \phi(x,\Omega,\tau) }{ S(r^{-1}(\tau) ) }Ê- \bar \psi(x,\Omega,r^{-1}(\tau)  d\Omega  \right) \frac{\delta \phi(x,\Omega',\tau)}{S (r^{-1}(\tau) )}  \frac{ dx d\Omega' d\tau }{r' (r^{-1}(\tau) )} =\\
 \int_{0}^{T_{R}} \int_{\X} \alpha_{1} \left( \int_{S^2}  \frac{ \phi(x,\Omega,\tau) }{ S(r^{-1}(\tau) ) }Ê- \bar \psi(x,\Omega,r^{-1}(\tau)) d\Omega \right) \delta \phi(x,\Omega',\tau) dx d\Omega' d\tau.
 \end{eqnarray*}
\par
\item 
Instead of applying the  transformation $r'(\epsilon) = \frac{1}{S(\epsilon)}$, we could also study the properties of an energy--dependent advection operator. If we introduce $\varphi = \psi f(\epsilon)$ with $f(\epsilon) = 1 / S(\epsilon)$ we obtain the family of operator 
$$ A(\epsilon) := f(\epsilon) \Omega\cdot\nabla_{x}.$$ 
Since $f>0$ this operator generates an evolution system $G(\epsilon,s)$  with domain 
of definition 
$$D(A(\epsilon))= \{ \psi \in L^2(\X\times\T):  - f(\epsilon) \Omega \nabla_{x} \psi \in L^2(\X\times\T), \psi|_{\Gamma_{-}} = 0 \}$$
and by 
\begin{equation*}
(G(\epsilon,s)\eta)(x,\Omega)\stackrel{def}{=}\eta(x- \Omega \cdot\int_{s}^{\epsilon}f(\varsigma) d\varsigma \,;\,\Omega) \; \forall  \eta \in C^0_{c}(\X). 
 \end{equation*}

\item 
In the beginning of Section \ref{sec:main result} we assumed that $\sigma_{t}$ and $\sigma_{s}$ are independent of the energy level $\epsilon.$ 
This is of course not a realistic assumption. Of course, formally, there will only be small changes to the  first--order optimality system (\ref{opti main}) in
the case of energy dependent coefficients. One simply replaces $\sigma_{t}$ and $\sigma_{s}$ through their energy dependent counterparts. However,
from an analytical point of view the semigroup theory presented here to solve the problem has to be extended in order to treat now  evolution equations. 
Furthermore, the transformation $ \epsilon \to r(\epsilon)$ used in order to establish the main result cannot be used to simplify the problem. All these points 
will be discussed in future work. 
 
\end{itemize}

\section*{Acknowledgments}
This work has been supported by DFG SPP1253 and KL 1105/14/2, DAAD D/08/11076 and RWTH Aachen Seed Funds 2008 and HE5386/6-1. 



\begin{thebibliography}{20}

\bibitem{HA} H. Amann, {\em Linear and Quasilinear Parabolic Problems}, Volume 1, Abstract Lienar Theory, 
Birkh\"auser Verlag, Basel, 1995

\bibitem{AydOliGod02}
E.~D. Aydin, C.~R.~E. Oliveira, and A.~J.~H. Goddard, \emph{A comparison
  between transport and diffusion calculations using finite element-spherical
  harmonics radiation transport method}, Med. Phys., Vol. {29}, 2002,
  2013--2023.

\bibitem{BelMai05}
N.~Bellomo and P.~K. Maini, \emph{Preface (special issue on cancer modelling)},
  Math. Mod. Math. Appl. Sci., Vol. {15}, No. 11, 2005, iii -- viii.

\bibitem{BelMai06}
N.~Bellomo and P.~K. Maini, \emph{Preface (special issue on cancer modelling)}, Math. Mod. Math.
  Appl. Sci. Vol. {16}, No. 7S, 2006, iii -- vii.

\bibitem{BelMai07}
N.~Bellomo and P.~K. Maini, \emph{Preface (special issue on cancer modelling)}, Math. Mod. Math.
  Appl. Sci. Vol. {17}, No. 11, 2007, iii -- vii.

\bibitem{BelLiMai08}
N.~Bellomo, N.~K.~Li, and P.~K. Maini, \emph{On the foundations of cancer modelling: selected topics, speculations, and perspectives}, Math. Mod. Math. Appl. Sci. Vol. {18}, No. 4, 2008, 593-646.

\bibitem{Bor98}
C.~B\"orgers, \emph{Complexity of {M}onte {C}arlo and deterministic
  dose-calculation methods}, Phys. Med. Biol., Vol. {43}, 1998, 517--528.

\bibitem{Bor99}
C.~B\"orgers, \emph{The radiation therapy planning problem}, {IMA} Volumes in
  Mathematics and its applications, vol. 110, Springer-Verlag, 1999.

\bibitem{BucBevRoa05}
K.~K. Bucci, A.~Bevan, and M.~Roach III, \emph{Advances in radiation therapy:
  conventional to 3d, to {IMRT}, to 4d, and beyond}, CA Cancer J. Clin.
  Vol. {55}, 2005, 117--134.


\bibitem{RL}
R.~Dautray and J.~L. Lions, \emph{Mathematical analysis and numerical methods
  for science and technology (v.6)}, Springer Verlag, Springer, 1993.


\bibitem{FHS} M. Frank, M. Herty and M. Sch\"afer. \,\textit{Optimal treatment plannig in radiotherapy based on Boltmann Transport Calculation}, Mathematical Models and Methods in Applied Sciences, Vol. 18(4), 2008, 573--592.


\bibitem{HenIzaSie06}
H.~Hensel, R.~Iza-Teran, and N.~Siedow, \emph{Deterministic model for dose
  calculation in photon radiotherapy}, Phys. Med. Biol. Vol. {51}, 2006,
  675--693.


\bibitem{HerPinSea07}
M.~Herty, R.~Pinnau, and M.~Seaid, \emph{Optimal control in radiative
  transfer}, Optimization Methods and Software, 2008.

\bibitem{GifHorWarFaiMou06}
K.~A. Gifford, J.~L.~Horton Jr., T.~A. Wareing, G.~Failla, and F.~Mourtada,
  \emph{Comparioson of a finite-element multigroup discrete-ordinates code with
  {M}onte {C}arlo for radiotherapy calculations}, Phys. Med. Biol. Vol. {51}
  2006, 2253--2265.

\bibitem{GM} Y. Giga Y. and T. Miyakawa, \textit{Solution in $L^r$ of Navier-Stokes initial value problem}, Arch.  Rat. Mech. Anal., Vol. 89, 1985, pp.~267-281.

\bibitem{KriSau05}
T. Krieger and O. Sauer, {\em {M}onte {C}arlo versus pencil-beam-/collapsed-cone-dose
calculation in a heterogeneous multi-layer phantom}, Phys. Med. Biol., Vol. 50,
2005, pp. 859--868





\bibitem{Lar97}
E.~W. Larsen, \emph{Tutorial: The nature of transport calculations used in
  radiation oncology}, Transp. theory Stat. Phys. Vol. 26, 1997, 739.

\bibitem{LarMifFraBru97}
E. W. Larsen and M. M. Miften and B. A. Fraass and I. A.
D. Bruinvis, {\em Electron dose calculations using the method of moments},
Med. Phys., Vol. 24, 1997, pp. 111--125




\bibitem{LN} L. Nirenberg, \textit{On elliptic partial differential equations}, Annali della Scoula Norm. Sup. Pisa, Vol. 13, 1959, pp.~115-162.



\bibitem{PA} A. Pazy, \textit{Semigroups of linear operators and applications to partial differential equations}. Applied Mathematical Sciences, Vol. 44,  Springer-Verlag,\,New York, Berlin,1983.


\bibitem{RR} M. Renardy and R.C. Rogers, {\em An Introduction to Partial Differential Equations}, Springer Texts in Applied Mathematics, Volume 13, Springer--Verlag, New York, 1996

\bibitem{SY} R. G. Sell and Y. You, {\em Dynamics of Evolutionary Equations}, Springer--Verlag, New York, Berlin, 2002.  


\bibitem{SheFerOliMac99}
D.~M. Shepard, M.~C. Ferris, G.~H. Olivera, and T.~R. Mackie, \emph{Optimizing
  the delivery of radiation therapy to cancer patients}, SIAM Rev. Vol. {41}
  (1999), 721--744.

%


\bibitem{HT1} H. Tanabe and T. Tanabe \textit{Functional Analysis for Partial Differential Parabolic Equations}, Marcel Dekker Publishers, 1996. 


\bibitem{TerKol02}
J.~Tervo and P.~Kolmonen, \emph{Inverse radiotherapy treatment planning model
  applying boltzmann-transport equation}, Math. Models. Methods. Appl. Sci.
  Vol. {12} (2002), 109--141.


\bibitem{TerKolVauHeiKai99}
J.~Tervo, P.~Kolmonen, M.~Vauhkonen, L.~M. Heikkinen, and J.~P. Kaipio, \emph{A
  finite-element model of electron transport in radiation therapy and related
  inverse problem}, Inv. Probl. Vol. {15} (1999), 1345--1361.

\bibitem{FT}
F.~Tr\"oltzsch, \emph{Optimale {S}teuerung partieller {D}ifferentialgleichungen
  - {T}heorie, {V}erfahren und {A}nwendungen}, Vieweg Verlag, 2005.

\bibitem{HT} H. Triebel, \textit{Interpolation Theory, Function Spaces, Differential Operators},  North--Holland Publishers,  1978.


\bibitem{FBW} F.B.  Weisler, \textit{Semilinear evolution equations in Banach spaces}, Journal of Functional Analysis, Vol. 32, 1979, pp.~ 277--296.




\end{thebibliography}
\end{document}